\documentclass[11pt,reqno]{amsart}
\usepackage{amssymb,latexsym}
\usepackage{graphicx}
\usepackage{fancyhdr,amssymb}
\usepackage{url}		

\calclayout

\makeatletter
\newcommand{\monthyear}[1]{%
  \def\@monthyear{\uppercase{#1}}}
\newcommand{\volnumber}[1]{%
  \def\@volnumber{\uppercase{#1}}}
\AtBeginDocument{%
\def\ps@plain{\ps@empty
  \def\@oddfoot{\@monthyear \hfil \thepage}%
  \def\@evenfoot{\thepage \hfil \@volnumber}}
\def\ps@firstpage{\ps@plain}
\def\ps@headings{\ps@empty
  \def\@evenhead{%
    \setTrue{runhead}%
    \def\thanks{\protect\thanks@warning}%
    \uppercase{Gaurav Bhatnagar}\hfil} 
  \def\@oddhead{%
    \setTrue{runhead}%
    \def\thanks{\protect\thanks@warning}%
    \hfill\uppercase{Analogues of a Fibonacci-Lucas identity}}%
  \let\@mkboth\markboth
  \def\@evenfoot{%
    \thepage \hfil \@volnumber}%
  \def\@oddfoot{%
    \@monthyear \hfil \thepage}%
  }%
\footskip=25pt
\pagestyle{headings}%
}
\makeatother

\theoremstyle{plain}
\numberwithin{equation}{section}
\newtheorem{thm}{Theorem}[section]
\newtheorem{theorem}[thm]{Theorem}

\newtheorem*{rems}{Remarks} 
\newenvironment{Remarks}{\begin{rems}\normalfont}{\end{rems}}

\makeatletter
\renewcommand{\@biblabel}[1]{#1.}
\makeatother

\newcommand{\qrfac}[2]{{\left({#1}; q \right)_{#2}}} 

\begin{document}

\monthyear{}
\volnumber{}

\title{Analogues of a Fibonacci-Lucas identity}

\author{Gaurav Bhatnagar
}
\thanks{Research supported in part by the Austrian Science Fund (FWF): F 5008-N15. 
}
 \thanks{Part of this work was done while the author was a Visiting Scientist at the Indian Statistical Institute, Delhi Centre, 
7 S.~J.~S.~Sansanwal Marg, Delhi 110016.}
\address{Fakult\"at f\"ur Mathematik,  Universit\"at Wien \\
Oskar-Morgenstern-Platz 1, 1090 Wien, Austria.}
\email{bhatnagarg@gmail.com}



\begin{abstract}
Sury's 2014 proof of an identity for Fibonacci and Lucas numbers (Identity 236 of Benjamin and Quinn's 2003 book: {\em Proofs that count: The art of combinatorial proof}) has excited a lot of comment. We give an alternate, telescoping, proof of this---and associated---identities and generalize them. We also give analogous identities for other sequences that satisfy a three-term recurrence relation. 
\end{abstract}

\maketitle 

\section{Lucas' 1876 identity and its associates}
Let $F_k$ and $L_k$ be the $k$th Fibonacci and $k$th Lucas numbers, respectively, for $k=0, 1, 2, \dots $. Both sequences of numbers satisfy a common recurrence relation, for $n=0, 1, 2, \dots$
$$x_{n+2}=x_{n+1}+x_n;  $$
with differing initial values: $F_0=0, F_1=1$ but $L_0=2, L_1=1$. 

Fibonacci numbers satisfy many identities. The classic 1876 identity is due to Lucas himself: 
$$F_0+F_1+F_2+\cdots +F_{n-1}=F_{n+1}-1,$$
which (in view of the relation $F_{k+1}+F_{k-1} = L_k$, for $k\geq 1$) can  be written as
\begin{equation}\label{lucas1876}
\sum_{k=0}^n \left(L_k - F_{k+1}\right) = F_{n+1},
\end{equation}
where we used $1=L_0-F_1$, $F_0=L_1-F_2$,  $F_1=L_2-F_3$, and so on. Next consider Identity 236 of Benjamin and Quinn~\cite[p.~131]{BQ2003}:
\begin{equation}\label{sury}
\sum_{k=0}^n 2^k L_k = 2^{n+1}F_{n+1}.
\end{equation}
A recent proof of this identity by Sury \cite{sury2014} 
excited a lot of comment.  Kwong \cite{kwong2014} gave an alternate proof, and Marques \cite{marques2015} (see also Martinjak \cite{martinjak2015}) found an analogous identity that replaces $2$ by $3$. Marques' identity can be written (in the form presented by Martinjak):
\begin{equation}\label{marques-martinjak}
\sum_{k=0}^n 3^k \left(L_k +F_{k+1}\right) = 3^{n+1}F_{n+1}.
\end{equation}
On the lines of \eqref{sury}, Martinjak gave an identity with alternating signs:
\begin{equation}\label{martinjak}
\sum_{k=0}^n \frac{(-1)^k}{2^k} L_{k+1}  = \frac{(-1)^n}{2^n}F_{n+1}.
\end{equation} 
Finally, the other classic identity (we are unaware of its provenance) which consists of alternating sums of Fibonacci numbers, see
 Koshy \cite[Identity 20, p.~88]{koshy2001}: 
$$F_2-F_3+F_4+\cdots +(-1)^{n+2}F_{n+2}=(-1)^{n}F_{n+1},$$
can be written as 
\begin{equation}\label{lucas-alternating}
\sum_{k=0}^n (-1)^{k}\left( L_{k+1}-F_k \right) = (-1)^{n}F_{n+1}.
\end{equation}

Now many of the classic Fibonacci identities follow from telescoping, including Lucas' classic identity \eqref{lucas1876}.
And all telescoping formulas can be derived from what we \cite{GB2011} called the Euler's Telescoping Lemma, see \eqref{eulerfinite2} below.
If we apply this approach to prove the above identities, we are quickly led to the following generalization, for $t$ an
indeterminate: 
\begin{equation}\label{gb-sury}
\sum_{k=0}^n t^k \left(L_k +(t-2)F_{k+1}\right) = t^{n+1}F_{n+1}.
\end{equation}
Observe that when $t=1, 2,$ and $3$, then \eqref{gb-sury} reduces to \eqref{lucas1876}, \eqref{sury} and \eqref{marques-martinjak}, respectively. What is remarkable is that \eqref{gb-sury} also contains as special cases the identities \eqref{martinjak} and  \eqref{lucas-alternating}. We can see this by replacing $t$ by $-1/2$, and  $t=-1$, respectively. But in this case one has to massage the sums  a bit using the relations between Fibonacci and Lucas numbers. 

Alternatively, one can directly prove, for $t\neq 0$:
\begin{equation}\label{gb-martinjak}
\sum_{k=0}^n \frac{(-1)^k}{t^k} \left(L_{k+1} +(t-2)F_{k}\right)  = \frac{(-1)^n}{t^n}F_{n+1}.
\end{equation} 
When $t=1$ and $2$, then \eqref{gb-martinjak} reduces to  \eqref{lucas-alternating} and 
 \eqref{martinjak} respectively. It is a nice exercise to show that \eqref{gb-sury} and \eqref{gb-martinjak} are equivalent. 

The purpose of this paper is to generalize the above identities to sequences which satisfy a recurrence relation of the form
$$x_{n+2}=a_nx_{n+1}+b_nx_n,$$
where $a_n$ and $b_n$ are sequences of indeterminates, which can be specialized to complex numbers or polynomials when required.   As examples, in addition to \eqref{gb-sury} and \eqref{gb-martinjak}, we will note many analogous identities for other sequences that are defined by this recurrence relation. See  \cite{GB2011} for more examples of this approach.  

\section{Euler's Telescoping Lemma and its application}
Euler's Telescoping Lemma can be  written as \cite[eq.~(2.2)]{GB2011}:
\begin{equation}\label{eulerfinite2}
\sum_{k=1}^n w_k \frac{u_1u_2\cdots u_{k-1}}{v_1v_2\cdots v_{k}}=
\frac{u_1u_2\cdots u_{n}}{v_1v_2\cdots v_{n}} -1,
\end{equation}
where $w_k=u_k-v_k$. Here the product $u_1\cdots u_{k-1}$ is considered to be equal to $1$ when $k=1$. 
 The $u_k$ and $v_k$ are sequences of indeterminates. Its proof is immediate. Replace $w_k$ by $u_k-v_k$ in the sum on the left hand side, and \eqref{eulerfinite2} follows by telescoping. I showed in \cite{GB2011} that an equivalent form of this identity characterizes telescoping sums. Thus all telescoping sums can be obtained as suitable special cases. 

For example, to obtain \eqref{gb-sury} set $u_k=tF_{k+1}$, $v_k=F_k$. Then
\begin{align*}
w_k &= u_k-v_k\\
&= tF_{k+1}-F_k\\
&=(F_{k+1}-F_k)+ (t-1)F_{k+1}\\
&=  F_{k-1}+ (t-1)F_{k+1}\\
&= (F_{k-1}+F_{k+1})+ (t-2)F_{k+1}\\
&=  L_{k}+(t-2)F_{k+1} .
\end{align*}
Now substitute for $u_k$, $v_k$, and $w_k$ in \eqref{eulerfinite2} to obtain \eqref{gb-sury}.

Similarly, to obtain \eqref{gb-martinjak}, set $u_k=F_{k+1}$, $v_k=-tF_k$ in \eqref{eulerfinite2} 
and observe that
$w_k= L_{k+1}+(t-2)F_k$. 

For sequences satisfying a more general three-term recurrence relation, we have the following generalization. 
\begin{theorem}\label{lucas-gen}
Let  $a_n$ and $b_n$ be sequences, with $a_n\neq 0$, $b_n\neq 0$, for all $n$. Let $t$ be an indeterminate.
Consider a sequence $x_n$ that satisfies (for $n\geq 0$):
$$x_{n+2}=a_nx_{n+1}+b_nx_{n}.$$
Then the following identities hold for $n=0,1,2,\dots$, provided the denominators are not $0$. 
\begin{align}
\sum_{k=1}^n   \frac{t^{k-1}}{a_0a_1\cdots a_{k-1}} 
 \frac{b_{k-1}x_{k-1}+(t-1)x_{k+1}}{x_1 } 
&=\frac{t^n}{a_0a_1\cdots a_{n-1}}  \frac{ x_{n+1}}{x_1}-1.
\label{sury-gen1}\\
\sum_{k=1}^n  \frac{(-1)^k}{t^k}  \frac{a_1a_2\cdots a_{k-1}}{b_1b_2\cdots b_k}  
 \frac{x_{k+2}+(t-1)b_kx_{k}}{x_1 } 
&=\frac{(-1)^n}{t^n} \frac{a_1a_2\cdots a_{n}}{b_1b_2\cdots b_n}  \frac{ x_{n+1}}{x_1}-1.
\label{martinjak-gen1}
\end{align}
\end{theorem}
\begin{Remarks} The special case $t=1$ of identity \eqref{martinjak-gen1} has appeared previously in \cite[eq.~(10.11)]{GB2011}. When $a_k=1=b_k$, then \eqref{sury-gen1} reduces to \eqref{gb-sury}, and \eqref{martinjak-gen1} reduces to \eqref{gb-martinjak}.
\end{Remarks}
\begin{proof} We derive \eqref{sury-gen1} from Euler's Telescoping Lemma \eqref{eulerfinite2}. Set 
$u_k =tx_{k+1},$  $v_k=a_{k-1}x_{k}.$
Note that $w_k=(b_{k-1}x_{k-1}+(t-1)x_{k+1})$. 
Substituting in \eqref{eulerfinite2}, we immediately obtain \eqref{sury-gen1}.

Next, set 
$u_k =a_kx_{k+1}$ and $v_k=-tb_kx_{k}$, so that 
$w_k=(x_{k+2}+(t-1)b_kx_{k})$. 
Substituting in \eqref{eulerfinite2}, we obtain \eqref{martinjak-gen1}.
\end{proof}

Next, we note down some examples of sequences that satisfy a three term recurrence relation. 

\section{Examples of analogous identities}
First consider the $k$th Pell number $P_k$ and the $k$th Pell-Lucas number $Q_k$. Both these sequences satisfy the recurrence relation, for $n=0, 1, 2, \dots$ 
$$x_{n+2}=2x_{n+1}+x_n,$$
and have initial values: $P_0=0$, $P_1=1$; and, $Q_0=2=Q_1$. 
 It is easy to check that $P_k$ and $Q_k$ satisfy the relation $$P_{k+1}+P_{k-1} = Q_{k}, \text{ for } k\geq 1. $$
Substituting $a_k=2$, $b_k=1$, and $x_k=P_k$ in \eqref{sury-gen1}, multiplying both sides by $t$, and using the above relation between the Pell and Pell-Lucas numbers, we obtain a Pell analogue of \eqref{gb-sury}:
\begin{equation*}
\sum_{k=0}^n \frac{t^k}{2^k} \left(Q_k +(t-2)P_{k+1}\right) =\frac{t^{n+1}}{2^n}P_{n+1}.
\end{equation*}
When  $t$ is replaced by $2t$, we have
\begin{equation}\label{pell-sury}
\sum_{k=0}^n t^k \left(Q_k +2(t-1)P_{k+1}\right) =2t^{n+1}P_{n+1}.
\end{equation}
Similarly, equation \eqref{martinjak-gen1} yields, after replacing $t$ by $2t$,
for $t\neq 0$:
\begin{equation}\label{pell-martinjak}
\sum_{k=0}^n  \frac{(-1)^k}{2 t^k} \left(Q_{k+1} +2(t-1)P_{k}\right)  = \frac{(-1)^n}{t^n}P_{n+1}.
\end{equation} 

Compare the $t=1$ case of \eqref{pell-sury}
\begin{equation}\label{pell-sury-t1}
\sum_{k=0}^n Q_k  =2 P_{n+1}.
\end{equation}
 with \eqref{lucas1876} and \eqref{sury}.
The identity
\begin{equation}\label{pell-martinjak-t1}
\sum_{k=0}^n  (-1)^k Q_{k+1}   = (-1)^n 2 P_{n+1},
\end{equation} 
obtained by setting $t=1$ in \eqref{pell-martinjak} is similarly very appealing. 

Next, we obtain identities for the Lucas numbers. Recall that Lucas numbers too satisfy the same recursion as the Fibonacci numbers, and have initial values $L_0=2$ and $L_1=1$. In this case, \eqref{sury-gen1} yields:
\begin{equation}\label{lucas-sury}
1+\sum_{k=0}^n {t^k} \left( L_{k-1} +(t-1)L_{k+1}\right) =t^{n+1}L_{n+1},
\end{equation}
where we have taken $L_{-1}=0$ in order to write the sum from $0$ to $n$. 
Similarly, equation \eqref{martinjak-gen1} yields,
for $t\neq 0$:
\begin{equation}\label{lucas-martinjak}
1+ \sum_{k=1}^n  \frac{(-1)^k}{t^k} \left(L_{k+2} +(t-1)L_{k}\right)  = \frac{(-1)^n}{t^n}L_{n+1}.
\end{equation} 

The Derangement Numbers $d_n$, counting the number of derangements---that is, the number of permutations on $n$ letters with no fixed points---are usually not considered to be analogous to the Fibonacci numbers. However, they do satisfy a three term recurrence relation. Since the formulas \eqref{sury-gen1} and \eqref{martinjak-gen1} involve division by $x_1$, and $d_1=0$, we have to make an adjustment. We consider the shifted derangement numbers $D_n$ which satisfy $D_n=d_{n+1}$. They are defined by: $D_0=0$, $D_1=1$ and for $n\geq 0$
$$D_{n+2}=(n+2)D_{n+1}+(n+2)D_n.$$
In this case $a_k=k+2=b_k$. Equation \eqref{sury-gen1} gives us
\begin{equation}\label{derangement-sury}
1+\sum_{k=1}^n \frac{t^{k-1}}{(k+1)!} \left((k+1)D_{k-1} +(t-1)D_{k+1}\right) =\frac{t^n}{(n+1)!}D_{n+1}.
\end{equation}
Similarly, equation \eqref{martinjak-gen1} yields,
for $t\neq 0$:
\begin{equation}\label{derangement-martinjak}
\sum_{k=0}^n  \frac{(-1)^k}{t^k} \left( \frac{D_{k+2} +(t-1)(k+2)D_{k}}{k+2}\right)  = \frac{(-1)^n}{t^n}D_{n+1}.
\end{equation} 
Identities for the derangement numbers are obtained by replacing $D_n$ by $d_{n+1}$.

Next, we consider the $q$-Fibonacci numbers considered first by Schur, and later studied by Andrews \cite{and1970, and1974, and2004}, Carlitz \cite{carlitz1974, carlitz1975}, and Smith \cite{smith2008}. We use the notation of Garrett \cite{garrett2004} who has studied these sequences combinatorially. The $q$-Fibonacci numbers are defined by $F_0^{(a)}(q)=0$, $F_1^{(a)}(q)=1$, and
$$F_{n+2}^{(a)}(q) = F_{n+1}^{(a)}(q) + q^{a+n} F_n^{(a)}(q).$$
In this case $a_k=1$, $b_k=q^{a+k}$ and \eqref{sury-gen1} yields
\begin{equation}\label{qfib-sury}
1+\sum_{k=1}^n t^{k-1} \left( q^{a+k-1}F^{(a)}_{k-1}(q) +(t-1)F^{(a)}_{k+1}(q)\right) =t^nF^{(a)}_{n+1}(q).
\end{equation}
Similarly, equation \eqref{martinjak-gen1} yields,
for $t\neq 0$:
\begin{align}
\sum_{k=0}^n  \frac{(-1)^k}{t^k} & q^{-ak}q^{-{k+1\choose 2}} \left( F^{(a)}_{k+2}(q) +(t-1)q^{a+k}F^{(a)}_{k}(q) \right) \notag
\\
&= \frac{(-1)^n}{t^n}q^{-na}q^{-{n+1\choose 2}}F^{(a)}_{n+1}(q). \label{qfib-martinjak}
\end{align}

The above identities may be regarded as $q$-analogs of \eqref{gb-sury} and \eqref{gb-martinjak}. As a final set of examples, we give another set of $q$-analogs of these two identities. 

We need the notation for $q$-rising factorials.
The $q$-rising factorial (for $q$ a complex number) is defined as the product:
$$\qrfac{a}{m} :=
\begin{cases}
1 &{\text{ if } m=0},\\
(1-a)(1-aq)\cdots (1-aq^{m-1}) &{\text{ if }} m\geq 1.\\
\end{cases}
$$

In Euler's Telescoping Lemma \eqref{eulerfinite2}, consider $u_k=(1-tq^{a+k-1})F^{(a)}_{k+1}(q)$ and 
$v_k=F^{(a)}_{k}(q)$. Then 
$$w_k=q^{a+k-1}\left( F^{(a)}_{k-1}(q)-tF^{(a)}_{k+1}(q)\right),$$
and we obtain:
\begin{equation}\label{q-sury}
1+\sum_{k=1}^n q^{a+k-1} \qrfac{tq^{a}}{k-1}
 \left( F^{(a)}_{k-1}(q) - t F^{(a)}_{k+1}(q)\right) =\qrfac{tq^{a}}{n} F^{(a)}_{n+1}(q).
\end{equation}
Next in \eqref{eulerfinite2}, 
set $u_k=F^{(a)}_{k+1}(q)$ and 
$v_k=(-1)q^{k+a}(1-tq^{-(a+k)})F^{(a)}_{k}(q)$. Then 
$$w_k= F^{(a)}_{k+2}(q)-tF^{(a)}_{k}(q),$$
and we obtain, after some simplification:
\begin{equation}\label{q-martinjak}
\sum_{k=0}^n \frac{1}{t^k} 
 \frac{ F^{(a)}_{k+2}(q) - t F^{(a)}_{k}(q)}{ \qrfac{t^{-1}q^{a+1}}{k}}=
\frac{1}{t^n} \frac{F^{(a)}_{n+1}(q)}{\qrfac{t^{-1}q^{a+1}}{n}}
\end{equation}

By now, the reader may well think that we have strayed quite far from the classic identities \eqref{lucas1876} and \eqref{lucas-alternating}.  However, given the uniform way they have been derived, one can see that all the identities presented here are related to each other. We emphasize that the above are only a small set of possible analogues. There are many sequences of numbers, and of polynomials, that satisfy such a three term recurrence. They all have analogues of Lucas' identity and its associates.

\section{Whats the point?} 
Any identity of the type 
$$\sum_{k=0}^n t_k = T_n$$
is a telescoping sum,  since one can write $T_k-T_{k-1}=t_k$ and then sum both sides. But how do you make it telescope? For many identities, we can use the WZ method 
(see Petkov{\u s}ek, Wilf and Zeilberger \cite{AeqB}) and let the computer find the telescoping.  For sequences defined by a three term recurrence relation, this may not always be possible, or practical. But you can use Euler's Telescoping Lemma to find the telescoping, without using a computer! This is the point of  \cite{GB2011}, which gives many examples of this approach, including several involving analogues of the Fibonacci numbers. 

So if you have found a Fibonacci type identity, you can use Euler's Telescoping Lemma to find analogous identities for other interesting sequences.

\subsubsection*{Note Added to proof} After submitting this paper, we found that T. Edgar (The Fibonacci Quarterly {\bf 54.1} (2016), 79) has independently given identity \eqref{gb-sury} for the case where $t$ is a natural number bigger than $1$.

\medskip

\noindent MSC2010: 11B39


\begin{thebibliography}{9999}

\bibitem{and1970} G. E. Andrews,  \emph{A polynomial identity which implies the Rogers--Ramanujan identities}, Scripta Math. \textbf{28} (1970), 297--305. 


\bibitem{and1974} G.~E.~Andrews, \emph{Combinatorial Analysis and Fibonacci Numbers}, Fibonacci Quart. 
\textbf{12 (2)}  (1974), 141--146.


\bibitem{and2004} G.~E.~Andrews, \emph{Fibonacci numbers and the Rogers--Ramanujan identities}, Fibonacci Quart. \textbf{42 (1)}  (2004), 3--19.

\bibitem{BQ2003} A.~T.~Benjamin and J.~J.~Quinn, \emph{Proofs that really count: The art of combinatorial proof}, {\it in} Dolciani Mathematical Exposition \# 27, Mathematical Association of America, Washington D.C.,  2003.

\bibitem{GB2011} G.~Bhatnagar, \emph{In praise of an elementary identity of Euler}, Electronic J. Combinatorics, \textbf{18 (2)}, \# P13 (2011), 44pp.


\bibitem{carlitz1974}  L.~Carlitz, \emph{Fibonacci notes. III. q-Fibonacci numbers}, Fibonacci Quart. 
\textbf{12 (4)} (1974), 317--322.

\bibitem{carlitz1975}  L.~Carlitz,  
\emph{Fibonacci notes. IV. q-Fibonacci polynomials}, Fibonacci Quart. \textbf{13 (2)}
 (1975), 97--102.


\bibitem{garrett2004} K.~C.~Garrett, \emph{Weighted Tilings and q-Fibonacci Numbers}, Preprint  (2004), 10pp.

\bibitem{koshy2001} T.~Koshy, \emph{Fibonacci and Lucas numbers with applications}, Wiley-Interscience, NY, 2001.

\bibitem{kwong2014} H.~Kwong, \emph{An alternate proof of Sury's Fibonacci-Lucas Relation}, Amer.~Math.~Monthly  \textbf{121 (6)} (2014), 514.

\bibitem{marques2015} D.~Marques, \emph{A new Fibonacci-Lucas Relation}, Amer.~Math.~Monthly  \textbf{122 (7)} (2015), 236.

\bibitem{martinjak2015} I.~Martinjak, \emph{Two extensions of Sury's identity}, arXiv:12080.01444v1 (2015), 2 pp.

\bibitem{smith2008} N.~O.~Smith, \emph{On an \lq uncounted' fibonacci identity and its q-analogue},
Fibonacci Quart. \textbf{46-47 (1)} (2008), 73--78.

\bibitem{AeqB} M.~Petkov{\u s}ek, H.~S.~Wilf, D.~Zeilberger: \emph{A=B},  A. K. Peters, Wellesley, MA, 1996.


\bibitem{sury2014} B.~Sury, \emph{A polynomial parent to a Fibonacci-Lucas relation}, Amer.~Math.~Monthly  \textbf{121 (3)} (2014), 236.


\end{thebibliography}
\end{document}